\documentstyle[11pt]{article}
\newenvironment{beweis}{{\it Proof.}\ }{\ $\ \ \ \Diamond$ \\\ }

 \newcounter{nsatz}[section]
 \newcounter{nlemma}[section]
 \newcounter{ndef}[section]
 \newcounter{nhyp}[section]
 \newcounter{nconjecture}[section]
 \newcounter{ncor}[section]
 \newcounter{nrem}[section]
 \newcounter{nexample}[section]
 \newcounter{nprop}[section]

 \newenvironment{nsatz}{\refstepcounter{nsatz}{\bf \arabic{section}.\arabic{nsatz}}\
                {\sc\bf Theorem.\ }\it}{\\\\ \rm}

 \newenvironment{nhyp}{\setcounter{nhyp}{\value{nsatz}}\refstepcounter{nhyp}
                \setcounter{nsatz}{\value{nhyp}}
                {\bf \arabic{section}.\arabic{nsatz}}\
                {\sc\bf Hypothesis.\ }}{\\\\ \rm}

 \newenvironment{ncor}{\setcounter{ncor}{\value{nsatz}}
                \refstepcounter{ncor}
                \setcounter{nsatz}{\value{ncor}}
                {\bf \arabic{section}.\arabic{nsatz}}\
                {\sc\bf Corollary.\ }\it}{\\\\ \rm}

 \newenvironment{nrem}{\setcounter{nrem}{\value{nsatz}}
                \refstepcounter{nrem}
                \setcounter{nsatz}{\value{nrem}}
                {\bf \arabic{section}.\arabic{nsatz}}\
                {\sc\bf Remark.\ }}{\\\\ \rm}

 \newenvironment{nprop}{\setcounter{nprop}{\value{nsatz}}
                \refstepcounter{nprop}
                \setcounter{nsatz}{\value{nprop}}
                {\bf \arabic{section}.\arabic{nsatz}}\
                {\sc\bf Proposition.\ }}{\\\\ \rm}

\begin{document}
\newcommand{\n}{{\mbox{\rm I$\!$N}}}
\newcommand{\z}{{\mbox{{\sf Z\hspace{-0.4em}Z}}}}
\newcommand{\R}{{\mbox{\rm I$\!$R}}}
\newcommand{\Q}{{\mbox{\rm I$\!\!\!$Q}}}
\newcommand{\C}{{\mbox{\rm I$\!\!\!$C}}}

\thispagestyle{empty}
\setlength{\parindent}{0pt}
\setlength{\parskip}{5pt plus 2pt minus 1pt}

\thispagestyle{empty}
\newcommand{\Irr}{{\mbox{\rm Irr}}}
\newcommand{\sIrr}{{\mbox{\scriptsize\rm Irr}}}
\newcommand{\Char}{{\mbox{\rm Char}}}
\mbox{\vspace{4cm}}
\vspace{4cm}
\begin{center}
{\bf \Large\bf  Counting characters in linear group actions\\}
\vspace{3cm}
by\\
\vspace{11pt}
Thomas Michael Keller\\
Department of Mathematics\\
Texas State University\\
601 University Drive\\
San Marcos, TX 78666\\
USA\\
e--mail: keller@txstate.edu\\
\vspace{1cm}
\end{center}
\thispagestyle{empty}
\newpage

\begin{center}
\parbox{12.5cm}{{\small
{\sc Abstract.}
Let $G$ be a finite group and $V$ be a finite $G$--module. We present upper bounds for the cardinalities of
certain subsets of $\Irr(GV)$, such as the set of those $\chi\in\Irr(GV)$ such that, for a fixed
$v\in V$, the restriction of $\chi$ to $\langle v\rangle$ is not a multiple of the regular
character of $\langle v\rangle$. These results might be useful in attacking the non--coprime
$k(GV)$--problem.\\
}}
\end{center}
\normalsize

\section{Introduction}\label{section0}
Let $G$ be a finite group and $V$ be a finite $G$--module of characteristic $p$.
If $(|G|,|V|)=1$, then in \cite[Theorem 2.2]{knoerr} R. Kn\"orr presented a beautiful argument
showing how to obtain strong upper bounds for $k(GV)$ (the number of conjugacy classes of $GV$) by using
only information on $C_G(v)$ for a fixed $v\in V$. Note that his result immediately implies the important
special case that if $G$ has a regular orbit on $V$ (i.e., there is a $v\in V$ with $C_G(v)=1$), then
$k(GV)\leq|V|$, which was a crucial result in the solution of the $k(GV)$--problem. In this note we give
a much shorter proof of this result (see \ref{prop2.1} below).\\
The main objective of the paper, however, is to modify and generalize Kn\"orr's argument in various directions
to include non--coprime situations. This way we obtain a number of bounds on certain subsets of $\Irr(GV)$,
such as the following:\\

{\bf Theorem A. }{\it Let $G$ be a finite group and let $V$ be a finite $G$--module of characteristic $p$. Let
$v\in V$ and $C=C_G(v)$ and suppose that $(|C|,|V|)=1$. Then the number of irreducible characters whose restriction to $\langle v\rangle$ is not a
multiple of the regular character of $\langle v\rangle$ is bounded above by
\[\sum_{i=1}^{k(C)}|C_V(c_i)|,\]
where the $c_i$ are representatives of the conjugacy classes of $C$.}\\

{\bf Theorem B. }{\it Let $G$ be a finite group and $V$ be a finite $G$--module. Let $g\in G$ be of prime
order not dividing $|V|$. Then the number of irreducible characters of $GV$ whose restriction to $A=\langle g\rangle$
is not a multiple of the regular character of $\langle g\rangle$ is bounded above by
\[|C_G(g)|\ n(C_G(g),C_V(g)),\]
where $n(C_G(g),C_V(g))$ denotes the number of orbits of $C_G(g)$ on $C_V(g)$.}\\

Stronger versions and refinements of these results are proved in the paper. It is hoped that these
results prove useful in solving the non--coprime $k(GV)$--problem, as discussed, for instance, in
\cite{kellerkgv} and \cite{guralnick-tiep}. Theorem A and B will be proved in Sections \ref{section2} and
\ref{section3} below
respectively. In Section \ref{section1}, we will generalize a recent result of P. Schmid
\cite[Theorem 2(a)]{schmid}
stating that in the
situation of the $k(GV)$--problem, if $G$ has a regular orbit on $V$, then $k(GV)=|V|$ can only hold if $G$
is abelian. We prove\\

{\bf Theorem C. }{\it Let $G$ be a finite group and $V$ a finite faithful $G$--module with
$(|G|,|V|)=1$. Suppose that $G$ has a regular orbit on $V$. Then
\[k(GV)\leq |V|-|G|+k(G).\]}
Our proof is different from the approach taken in \cite{schmid}, and we actually will prove a slightly
stronger result including some non--coprime actions.\\

Notation: If the group $A$ acts on the set $B$, we write $n(A,B)$ for the number of orbits of $A$ on $B$.
All other notation is standard or explained along the way.\\

\section{${\bf k(GV)=|V|}$ and regular orbits}\label{section1}

In this paper we often work under the hypothesis of the $k(GV)$--problem which is the following.\\

\begin{nhyp}\label{hyp1.1}
Let $G$ be a finite group and let $V$ be a finite faithful $G$--module such that $(|G|,|V|)=1$.
Write $p$ for the characteristic of $V$.
\end{nhyp}

In \cite[Theorem 2(a)]{schmid} P. Schmid proved that under \ref{hyp1.1}, if $G$ has a regular
orbit on $V$, $V$ is irreducible, and $k(GV)=|V|$, then $G$ is abelian, and from this it follows
easily that either $|G|=1$ and $|V|=p$, or $G$ is cyclic of order $|V|-1$. The proof in \cite{schmid}
is somewhat technical.\\
The goal of this section is to give a short proof of a generalization of Schmid's result based
on a beautiful argument of Kn\"orr \cite{knoerr}. We word it in such a way that we even do not need
the coprime hypothesis, so that the result may even be useful to study the non--coprime $k(GV)$--problem.
To do this, for any group $X$ and $x\in X$ we introduce the set
\[\Irr(X,x)=\{\chi\in\Irr(G)|\ \chi|_{\langle x\rangle} \mbox{ is not an integer multiple of the regular
character of }\langle x\rangle\}\]
and write
\[k(X,x)=|\Irr(X,x)|.\]

\begin{nsatz}\label{satz1.2}
Let $G$ be a finite group and let $V$ be a finite $G$--module such
that $G$ possesses a regular orbit on $V$. Let $v\in V$ be a representative of such an orbit. Then
\[k(GV,v)\leq |V|-|G|+k(G)\]
\end{nsatz}
\begin{beweis}
Let $p$ be the characteristic if $V$.
We proceed exactly as in Case (ii) of the proof of \cite[Theorem 2.2]{knoerr}.
Write $C=C_G(v)$. As $C=1$, we see that for $A=\langle v\rangle$ we trivially have that $|C|$ and $|A|$ are coprime,
and so that proof yields
\[(1) \quad (p-1)|V|=\sum_{\tau\in\sIrr(GV)}(\tau\eta,\tau)_A\]
where $\eta$ is the character of $A$ defined by $\eta=p1_A-\rho_A$ with $\rho_A$ being the regular
character of $A$. Now for any $\tau\in\Irr(GV)$ we have
\begin{eqnarray*}
(2)\quad (\tau\eta,\tau)_A  &=&\frac{1}{|A|}\sum_{a\in A}\tau(a)(p-\rho_A(a))\overline{\tau(a)}\\
                            &=&\sum_{1\not=a\in A}|\tau(a)|^2\left\{\begin{array}{ll}
                                                    =0 &\mbox{if }\tau|_A\mbox{ is an integer multiple of }\rho_A\\
                                                    \geq p-1 &\mbox{otherwise}
                                                                 \end{array}\right.
\end{eqnarray*}
where the last step follows from \cite[Corollary 4]{robinson}.
Next observe that if $\tau\in\Irr(GV)$ with $V\leq\ker\tau$, then $\tau\in\Irr(G)$ and clearly
$\tau|_A$ is not a multiple of $\rho_A$, and then clearly
\[(3)\quad (\tau\eta,\tau)_A=\sum_{1\not=a\in A}|\tau(a)|^2=\sum_{1\not=a\in A}\tau(1)^2=(p-1)\tau(1)^2.\]
Thus with (1), (2), and (3) we get
\begin{eqnarray*}
(p-1)|V| &=&\sum_{\tau\in\sIrr(G)}(\tau\eta,\tau)_A+\sum_{{\scriptsize\begin{array}{l}
                                                        \tau\in\sIrr(GV),\\
                                                        V\not\leq\ker\tau
                                                        \end{array} }}
                                                  (\tau\eta,\tau)_A\\
         &\geq&\sum_{\tau\in\sIrr(G)}(p-1)\tau(1)^2+(k(GV,v)-k(G))(p-1)
\end{eqnarray*}
which yields
\[|V|\geq\sum_{\tau\in\sIrr(G)}\tau(1)^2+k(GV,v)-k(G) =|G|+k(GV,v)-k(G).\]
This implies the assertion of the theorem, and we are done.
\end{beweis}

The following consequence implies Schmid's result \cite[Theorem 2(a)]{schmid}.\\

\begin{ncor}\label{cor1.3}
Assume \ref{hyp1.1} and that $G$ has a regular orbit on $V$. Then
\[k(GV)\leq |V|-|G|+k(G).\]
In particular, if $k(GV)=|V|$, then $G$ is abelian.
\end{ncor}
\begin{beweis}
By Ito's theorem and as $(|G|,|V|)=1$, we know that $\chi(1)$ divides $|G|$ for every $\chi\in\Irr(GV)$,
so in particular $p$ does not divide $\chi(1)$. Thus for any $v\in V^\#$ we see that
$\chi|_{\langle v\rangle}$ cannot be an integer multiple of $\rho_{\langle v\rangle}$. Therefore
$k(GV,v)=k(GV)$. Now the assertion follows from \ref{satz1.2}.
\end{beweis}

\section{Bounds for ${\bf k(GV)}$}\label{section2}

In this section we study more variations of Kn\"orr's argument in \cite[Theorem 2.2]{knoerr} and
generalize it to some non-coprime situations.\\
We begin, however, by looking at a classical application of it. An important and immediate
consequence of Kn\"orr's result is that if under \ref{hyp1.1} $G$ has a regular orbit on $V$, then
$k(GV)\leq |V|$. This important result can be obtained in the following shorter way.\\

\begin{nprop}\label{prop2.1}
Let $G$ be a finite group and let $V$ be a finite faithful $G$--module. Let
$v\in V$. Then
\[k(GV,v)\ \leq\ |C_G(v)| |V|,\]
in particular, if $(|G|,|V|)=1$ and $G$ has a regular orbit on $V$, then $k(GV)\leq |V|$.
\end{nprop}
\begin{beweis}
Put $A=\langle v\rangle$. If $\tau\in\Irr(GV,v)$, then by \cite[Corollary 4]{robinson}
we know that $\sum\limits_{1\not= a\in A} |\tau(a)|^2\geq p-1$. With this and well--known
character theory we get
\begin{eqnarray*}
(p-1)k(GV,v)&\leq&k(GV,v)\min\limits_{\tau\in\sIrr(GV,v)}\left(
                  \sum\limits_{1\not=a\in A}|\tau(a)|^2\right)\\
            &\leq&\sum\limits_{\tau\in\sIrr(GV)}\sum\limits_{1\not=a\in A}|\tau(a)|^2\\
            &=&   \sum\limits_{1\not=a\in A}\sum\limits_{\tau\in\sIrr(GV)}\tau(a)\overline{\tau(a)}\\
            &=&   \sum\limits_{1\not=a\in A}|C_{GV}(a)|\\
            &=&   \sum\limits_{1\not=a\in A}|C_G(v)| |V|\\
            &=&   (p-1)|C_G(v)| |V|
\end{eqnarray*}
This implies the first result. If $(|G|,|V|)=1$, then by Ito's result $\tau(1)\big| |G|$
for all $\tau\in\Irr(GV)$, so $p$ cannot divide $\tau(1)$, and thus $k(GV,v)=k(GV)$,
and the second result now follows by choosing $v$ to be in a regular orbit of $G$ on $V$.
\end{beweis}

Now we turn to generalizing Kn\"orr's argument. We discuss various ways to do so.\\

\begin{nrem}\label{rem2.1}
Let $G$ be a finite group and let $V$ be a finite faithful $G$--module of characteristic $p$.
Let $v\in V$ and put $C=C_G(v)$ and $A=\langle v\rangle$.
Let
\[\Irr(GV,C,v):=\Irr_0(GV):=\Irr(GV) - \{\chi\in\Irr(GV)\ |\ \chi|_{C\times\langle v\rangle}=
\tau\times\rho_A\mbox{ for a character $\tau$ of $C$}\}\]
and
\[\Irr_{p'}(GV)=\{\chi\in\Irr(GV)\ |\ p\mbox{ does not divide }\chi(1)\},\]
so that clearly $\Irr_{p'}(GV)\subseteq\Irr_0(GV)$.\\
Note that if $(|G|,|V|)=1$, then by Ito $\Irr(GV)=\Irr_{p'}(GV)$.
\end{nrem}

To work towards our next result, we again proceed somewhat similarly as in \cite[Theorem 2.2]{knoerr}.
In the following we work under the hypothesis that $(|C|,|V|)=1$.
Let $N=N_G(A)$.
Then $|N:C|$ divides $p-1$. Moreover, from Kn\"orr's proof we know that if $c_i$ ($i=1,\ldots,k(C)$)
with $c_1=1$ are representatives of the conjugacy classes of $C$ and $a_j$ ($j=1,\ldots,\frac{p-1}{|N:C|}$)
are representatives of the $N$--conjugacy classes of $A-1$ then, the $c_ia_j$ are representatives
of those conjugacy classes of $GV$ which intersect $C\times (A-1)$ nontrivially.\\
Moreover recall from Kn\"orr's proof that for $c\in C$, $1\not= a\in A$,
$g\in G$, $u\in V$ we know that
\[(ca)^{gu}\in C\times A\mbox{ if and only if }g\in N\mbox{ and }u\in C_V(c^g).\]
Now define a character $\eta$ on $C\times A$ by $\eta=1_C\times(p1_A-\rho_A)$.\\
Then for $c\in C$, $a\in A$ we have
\[\eta(ca)\ =\ \left\{\begin{array}{ll}
p,&\mbox{if }a\not=1\\
0,&\mbox{if }a=1
\end{array}\right.\]
Therefore $\eta^{GV}$ vanishes on all conjugacy classes of $GV$ which intersect $C\times(A-1)$
trivially, whereas for $c\in C$, $1\not=a\in A$ we have that
\begin{eqnarray*}
\eta^{GV}(ca)&=&\frac{1}{|C\times A|}\sum_{{\scriptsize\begin{array}{l}g\in G\\
                                                           u\in V\end{array}}}
                                        \dot{\eta}((ca)^{gu})\\
             &=&\frac{1}{p|C|}\sum_{g\in N}\sum_{u\in C_V(c^g)}\eta(c^ga^g)\\
             &=&\frac{1}{p|C|}\sum_{g\in N}|C_V(c^g)|p\\
             &=&\frac{1}{|C|}\sum_{g\in N}|C_V(c)|\\
             &=&|N:C|\ |C_V(c)|.
\end{eqnarray*}
Thus if $x_i$ ($i=1,\ldots,k(GV)$) are representatives of the conjugacy classes of $GV$, then we get
\begin{eqnarray*}
\sum_{i=1}^{k(GV)}\eta^{GV}(x_i)&=&\sum_{i=1}^{k(C)}\sum_{j=1}^{\frac{p-1}{|N:C|}}\eta^{GV}(c_ia_j)\\
                                &=&\sum_{i=1}^{k(C)}\sum_{j=1}^{\frac{p-1}{|N:C|}}|N:C|\ |C_V(c_i)|\\
                                &=&\frac{p-1}{|N:C|}|N:C|\sum_{i=1}^{k(C)}|C_V(c_i)|\\
                                &=&(p-1)\sum_{i=1}^{k(C)}|C_V(c_i)|,
\end{eqnarray*}
and thus
\begin{eqnarray*}
(p-1)\sum_{i=1}^{k(C)}|C_V(c_i)|&=&\sum_{i=1}^{k(GV)}\eta^{GV}(x_i)\ =\
                                       \sum_{\tau\in\sIrr(GV)}(\tau\eta^{GV},\tau)_{GV}\\
                                &=&\sum_{\tau\in\sIrr(GV)}(\tau\eta,\tau)_{C\times A}\quad (1).
\end{eqnarray*}

Now if $\tau\in\Irr(GV)$, as in \cite{knoerr} write
\[\tau|_{C\times A}\ =\ \sum_{\lambda\in\sIrr(A)}\tau_\lambda \times\lambda\quad (2)\]
where $\tau_\lambda$ is a character of $C$ or $\tau_\lambda =0$.\\
Then as in \cite{knoerr} we see that
\begin{eqnarray*}
(\tau\eta,\tau)_{C\times A}&=&\frac{1}{|C\times A|}\sum_{{\scriptsize\begin{array}{l}c\in C\\
                                                                         a\in A
                                                                         \end{array}}}
                                                      \tau(ca)\eta(ca)\overline{\tau(ca)}\\
                           &=&\frac{1}{|C|}\sum_{{\scriptsize\begin{array}{l}c\in C\\
                                                         1\not=a\in A
                                                         \end{array}}}
                                                \tau(ca)\overline{\tau(ca)}\\
                           &=&\sum_{\lambda<\mu}((\tau_\lambda-\tau_\mu),(\tau_\lambda-\tau_\mu))_C\quad (3)
\end{eqnarray*}
where "$\leq$" is some arbitrary ordering on $\Irr(A)$.\\
Now if $\tau_\lambda-\tau_\mu$ is a nonzero multiple of $\rho_C$, then
\[(\tau_\lambda-\tau_\mu,\tau_\lambda-\tau_\mu)_C\ \geq\ |C|\quad (4)\]
and thus
\[(\tau\eta,\tau)_{C\times A}\ \geq\ |C|.\]
Moreover, note that if $\tau\in\Irr_0(GV)$, then not all $\tau_\lambda-\tau_\mu$ can be equal to 0 as otherwise
from (2) we see that $\tau|_{C\times A}$ would be equal to $\tau_\lambda\times \rho_A$ for any $\lambda$.
So we can partition the set $\Irr (A)$ into two disjoint nonempty subsets $\Lambda_1=\{\lambda\in\Irr(A)\ |\
\tau_\lambda=\tau_1\}$ and $\Lambda_2=\{\lambda\in\Irr(A)\ |\ \tau_\lambda\not=\tau_1\}$, and thus as in
\cite{knoerr} we see that $|\Lambda_1|\ |\Lambda_2|\geq p-1$, so there are at least $p-1$ pairs
$\lambda,\mu\in\Irr(A)$ such that $\tau_\lambda-\tau_\mu\not=0$. Thus
\[(\tau\eta,\tau)_{C\times A}\geq p-1\mbox{ for all }\tau\in\Irr_{0}(GV).\quad (5)\]
Therefore by (1) and (5) we get that
\[(p-1)\sum_{i=1}^{k(C)}|C_V(c_i)|=\sum_{\tau\in\sIrr(GV)}(\tau\eta,\tau)_{C\times A}\geq
\sum_{\tau\in\sIrr_0(GV)}(\tau\eta,\tau)_{C\times A}\geq(p-1)|\Irr_0(GV)|\]
and thus
\[|\Irr_0(GV)|\leq\sum_{i=1}^{k(C)}|C_V(c_i)|.\quad (6)\]
From now on we assume that $C>1$.\\
Now we repeat the arguments of this proof, but replace $\eta$ by
\[\eta_1\ =\ (|C| 1_C-\rho_C)\times(p1_A-\rho_A),\]
so for $c\in C$ and $a\in A$ we have
\[\eta_1(ca)\ =\ \left\{\begin{array}{ll}
                       |C|p&\mbox{if }c\not=1\mbox{ and }a\not=1\\
                       0   &\mbox{if }c=1\mbox{ or }a=1
                       \end{array}\right.\]
Now from the above we know that the $c_ia_j$
($i=2,\ldots,k(C)$, $j=1,\ldots,\frac{p-1}{|N:C|}$) are representatives of those conjugacy classes which
intersect $(C-1)\times (A-1)$ nontrivially.\\
Clearly $\eta_1^{GV}$ vanishes on all conjugacy classes of $GV$ which intersect $(C-1)\times (A-1)$
trivially, whereas for $1\not=c\in C$, $1\not=a\in A$, if $(|C|,|V|)=1$, we have that
\begin{eqnarray*}
\eta_1^{GV}(ca)&=&\frac{1}{|C\times A|}\sum_{{\scriptsize\begin{array}{l}g\in G\\
                                                            u\in V
                                                            \end{array}}}\dot{\eta}_1((ca)^{gu})\\
              &=&\frac{1}{p|C|}\sum_{g\in N}\sum_{u\in C_V(c^g)}\eta_1(c^ga^g)\\
              &=&|N|\ |C_V(c)|.
\end{eqnarray*}
Next we conclude that
\[\sum_{i=1}^{k(GV)}\eta_1^{GV}(x_i)=\sum_{i=2}^{k(C)}\sum_{j=1}^{\frac{p-1}{|N:C|}}\eta_1^{GV}(c_ia_j)
=(p-1)|C|\sum_{i=2}^{k(C)}|C_V(c_i)|,\]
and so as in (1) we see that
\[(p-1)|C|\sum_{i=2}^{k(C)}|C_V(c_i)|=\sum_{\tau\in\sIrr(GV)}(\tau\eta_1,\tau)_{C\times A}\quad (7).\]
Now with (2) similarly as in \cite{knoerr} we see that
\begin{eqnarray*}
(\tau\eta_1,\tau)_{C\times A}&=&\frac{1}{|C\times A|}\sum_{{\scriptsize\begin{array}{l}c\in C\\
                                                                           a\in A
                                                                           \end{array}}}
                                                          \tau(ca)\eta_1(ca)\overline{\tau(ca)}\\
                            &=&\sum_{{\scriptsize\begin{array}{l}1\not=c\in C\\
                                                     1\not=a\in A
                                                     \end{array}}}
                                   \tau(ca)\overline{\tau(ca)}\\
                            &=&\sum_{{\scriptsize\begin{array}{l}1\not=c\in C\\
                                                     1\not=a\in A
                                                     \end{array}}}\sum_{\lambda\in\sIrr (A)}
                                    \tau_\lambda(c)\lambda(a)\sum_{\mu\in\sIrr (A)}
                                    \overline{\tau_\mu(c)}\overline{\mu(a)}\\
                            &=&\sum_{\lambda,\mu\in\sIrr (A)}\sum_{1\not=c\in C}
                                   \tau_\lambda(c)\overline{\tau_\mu(c)}\sum_{1\not=a\in
                                   A}\lambda(a)\overline{\mu(a)}\\
                            &=&(p-1)\sum_{\lambda\in\sIrr (A)}\sum_{1\not=c\in C}\tau_\lambda(c)\overline{\tau_\lambda(c)}-
                               \sum_{{\scriptsize\begin{array}{l}\lambda,\mu\in \sIrr (A)\\
                                                     \lambda\not=\mu
                                     \end{array}}}
                               \sum_{1\not=c\in C}\tau_\lambda(c)\overline{\tau_\mu(c)}\\
                            &=&p\sum_{\lambda\in \sIrr (A)}\sum_{1\not=c\in C}\tau_\lambda(c)\overline{
                               \tau_\lambda(c)}-\sum_{\lambda,\mu\in\sIrr (A)}\sum_{1\not=c\in C}
                               \tau_\lambda(c)\overline{\tau_\mu(c)}\\
                            &=&\sum_{\lambda<\mu}\sum_{1\not=c\in C}(\tau_\lambda(c)-\tau_\mu(c))
                               (\overline{\tau_\lambda(c)}-\overline{\tau_\mu(c)})\\
                            &=&\sum_{\lambda<\mu}\sum_{1\not=c\in C}|\tau_\lambda(c)-\tau_\mu(c)|^2\quad (8)
\end{eqnarray*}
for some arbitrary ordering $\leq$ on $\Irr (A)$.\\
Now recall that if $\tau\in\Irr_0(GV)$, then not all of the $\tau_\lambda-\tau_\mu$ can be 0. So choose
$\lambda,\mu\in\Irr(C)$ such that $\tau_\lambda-\tau_\mu\not=0$. If all the $\tau_\mu$ ($\mu\in\Irr(A)$)
are integer multiples of $\rho_C$ then put $\Lambda_1=\{\phi\in\Irr(A)\ |\ \tau_\phi=\tau_\lambda\}$
and $\Lambda_2=\{\phi\in\Irr(A)\ |\ \tau_\phi\not=\tau_\lambda\}$, so $\Lambda_1\not=\emptyset$ and
$\Lambda_2\not=\emptyset$
and from $0\leq(|\Lambda_1|-1)(|\Lambda_2|-1)$ we clearly deduce that $|\Lambda_1||\Lambda_2|\geq p-1$,
so there are at least $p-1$ pairs $(\phi_1,\phi_2)\in\Irr(A)\times\Irr(A)$ such that
$\tau_{\phi_1}-\tau_{\phi_2}$ is a nonzero multiple of $\rho_C$.\\
So next we assume that $\tau_\lambda$ is not a multiple of $\rho_C$.\\
Then put
\[\Gamma_1=\{\phi\in\Irr(A)\ |\ \tau_\lambda-\tau_\phi\mbox{ is a multiple of }\rho_C\}\]
and
\[\Gamma_2=\{\phi\in\Irr(A)\ |\ \tau_\lambda-\tau_\phi\mbox{ is not a multiple of }\rho_C\}.\]
Clearly $\lambda\in\Gamma_1$, so $\Gamma_1\not=\emptyset$. If $\Gamma_2=\emptyset$, then
$\Irr(A)=\Gamma_1$, and if we define $\Lambda_1$, $\Lambda_2$ as in the previous argument, we see that
there are at least $(p-1)$ pairs $(\phi_1,\phi_2)\in\Irr (A)\times\Irr (A)$ such that
$\tau_{\phi_1}-\tau_{\phi_2}$ is a nonzero multiple of $\rho_C$.\\
So now suppose $\Gamma_2\not=\emptyset$. Then $|\Gamma_1|+|\Gamma_2|=p$, and if $\phi_1\in\Gamma_1$ and
$\phi_2\in\Gamma_2$, then $\tau_{\phi_1}-\tau_{\phi_2}=(\tau_{\phi_1}-\tau_\lambda)+(\tau_\lambda-\tau_{\phi_2})$
clearly is not a multiple of $\rho_C$, and by the same argument as used before we see that
$|\Gamma_1||\Gamma_2|\geq p-1$, so there are at least $(p-1)$ pairs $(\phi_1,\phi_2)\in\Irr (A)\times\Irr (A)$
such that $\tau_{\phi_1}-\tau_{\phi_2}$ is not a multiple of $\rho_C$.\\
Altogether we thus have shown that for any $\tau\in\Irr_0(GV)$ one of the following holds:\\
(A) There are at least $(p-1)$ pairs $(\phi_1,\phi_2)\in\Irr (A)\times\Irr (A)$ such that
\[\tau_{\phi_1}-\tau_{\phi_2}\mbox{ is a nonzero multiple of }\rho_C,\mbox{ or}\]
(B) there are at least $(p-1)$ pairs $(\phi_1,\phi_2)\in\Irr (A)\times\Irr (A)$ such that
\[\tau_{\phi_1}-\tau_{\phi_2}\mbox{ is not a multiple of }\rho_C.\]
Now it remains to consider two cases:\\
Case 1: At least half of the $\tau\in\Irr_0(GV)$ satisfy (A).\\

Then for any of these $\tau$ by (3) and (4) we have
\[(\tau\eta,\tau)_{C\times A}=\sum_{\lambda<\mu}((\tau_\lambda-\tau_\mu),(\tau_\lambda-\tau_\mu))_C\geq
(p-1)|C|\]
and so by (1) we see that
\[(p-1)\sum_{i=1}^{k(C)}|C_V(c_i)|\geq\sum_{\tau\in\sIrr_0(GV)}(\tau\eta,\tau)_{C\times A}\geq
\frac{1}{2}|\Irr_0(GV)|(p-1)|C|\]
which implies
\[|\Irr_0(GV)|\leq\frac{2}{|C|}\sum_{\i=1}^{k(C)}|C_V(c_i)|\quad (9).\]
Case 2: At least half of the $\tau\in\Irr_0(GV)$ satisfy (B).\\

Then for any of these $\tau$ by (8) and \cite[Corollary 4]{robinson} we have
\[(\tau\eta_1,\tau)_{C\times A}\geq (p-1)(k(C)-1).\]
Thus by (7) we have that
\[(p-1)|C|\sum_{i=2}^{k(C)}|C_V(c_i)|\geq\sum_{\tau\in\sIrr_0(GV)}(\tau\eta_1,\tau)_{C\times A}\geq
\frac{1}{2}|\Irr_0(GV)|(p-1)\cdot (k(C)-1)\]
whence
\[|\Irr_0(GV)|\leq\frac{2|C|}{k(C)-1}\sum_{i=2}^{k(C)}|C_V(c_i)|\quad (10).\]
Now we drop the assumption $(|C|,|V|)=1$ and work towards a general bound for $|\Irr_0(GV)|$.\\
For this, fix $g_0\in $C such that $g_0$ is of prime order $q$ and put $C_0=\langle g_0\rangle$
and $N_0=N_G(C_0)$. Trivially there are at most $|C_0|(p-1)=q(p-1)$ conjugacy classes of $GV$ that
intersect $C_0\times (A-1)$ nontrivially, and given $1\not=c\in C_0$, $1\not=a\in A$, we see that
for $g\in G$, $u\in V$
\[(ca)^{gu}=c^g[c^g,u]a^g\ \in\ C_0\times A\mbox{ first implies }c^g\in C_0,\mbox{ i.e., }g\in N_0,\]
and for each fixed $g\in N_0$, the equation $[c^g,u]a^g\in A$ implies
$[c^g,u]\in Aa^{-g}$ which has at most $|C_V(c^g)|\ |Ag^{-1}|=p|C_V(g_0)|$ solutions $u$.\\
Moreover, if $c=1$, then
\[(ca)^{gu}=a^{gu}=a^g\mbox{ implies }g\in N_G(A)=N\mbox{ and }u\in V.\]
Now we define the character $\eta_2$ on $C_0\times A$ by $\eta_2=1_{C_0}\times(p1_A-\rho_A)$.
Thus $\eta_2^{GV}$ vanishes on all conjugacy classes of $GV$ which intersect $C_0\times(A-1)$
trivially, whereas for $1\not=c\in C_0$, $1\not=a\in A$ we get
\begin{eqnarray*}
\eta_2^{GV}(ca)&=&\frac{1}{|C_0\times A|}\sum_{{\scriptsize\begin{array}{l}
                                                  g\in G\\
                                                  u\in V
                                                  \end{array}}}
                                         \dot{\eta}\left( (ca)^{gu}\right)\\
           &\leq&\frac{1}{qp}\sum_{g\in N_0}p|C_V(g_0)|p\\
           &=&   \frac{p}{q}|N_0||C_V(g_0)|,
\end{eqnarray*}
and for $c=1$, $1\not=a\in A$ we get
\[\eta_2^{GV}(ca)=\eta_2^{GV}(a)=\frac{1}{qp}\sum_{q\in N}|V|p=\frac{1}{q}|N||V|.\]
Thus if $x_i$ ($i=1,\ldots,k(GV)$) are representatives of the conjugacy classes of $GV$, then
\[\sum_{i=1}^{k(GV)}\eta_2^{GV}(x_i)\leq(p-1)\frac{1}{q}|N||V|+(q-1)(p-1)\frac{p}{q}|N_0||C_V(g_0)|\]
and as in (1) we see that
\[\sum_{i=1}^{k(GV)}\eta_2^{GV}(x_i)=\sum_{\tau\in\sIrr(GV)}(\tau\eta_2,\tau)_{C_0\times A}.\]
Now arguing as in (2), (3), (5) and (6) above will yield
\[|\Irr_{p'}(GV)|\leq k(GV,v)\leq |\Irr(GV,C_0,v)|\leq\frac{1}{q}(|N||V|+(q-1)p|N_0||C_V(g_0)|),\]
where $\Irr(GV,C_0,v)$ is as defined at the beginning of \ref{rem2.1}.
Putting the main results together, altogether we have proved the following:\\

\begin{nsatz}\label{satz2.2}
Let $G$ be a finite group and let $V$ be a finite faithful $G$--module of characteristic $p$.
Let $v\in V$ and put $C=C_G(v)$. If $c_i$ ($i=1,\ldots,k(C)$) are representatives of the
conjugacy classes of $C$, then the following hold:\\

(a) If $(|C|,|V|)=1$, then
\[|\Irr_0(GV)|\leq\sum_{i=1}^{k(C)}|C_V(c_i)|\]
and if $C>1$, then
\[|\Irr_0(GV)|\leq\max\left\{\frac{2}{|C|}\sum_{i=1}^{k(C)}|C_V(c_i)|,\frac{2|C|}{k(C)-1}\sum_{i=2}^{k(C)}
|C_V(c_i)|\right\}\]
(b) If $(|G|,|V|)=1$, then
\[\Irr_0(GV)=\Irr(G), \mbox{ so }k(GV)=|\Irr_0(GV)|\]
and the bounds in (a) hold true for $k(GV)$ instead of $|\Irr_0(GV)|$.\\

(c) In general, if $g\in C$ such that $o(g)=q$ is a prime, then
\[|\Irr_{p'}(GV)|\leq k(GV,v)\leq
\frac{1}{q}\Big(|N_G(\langle v\rangle )||V|+(q-1)p|N_G(\langle g\rangle)||C_V(g)|\Big).\]
\end{nsatz}

\section{The dual approach}\label{section3}

In the previous section, we always fixed $v\in V$ and obtained bounds on the size of suitable subsets
of $\Irr(GV)$ in terms of properties of the action of $C_G(v)$ on $V$. In this section we consider a
"dual" approach:\\
We fix $g\in G$ and find bounds in terms of the action of $C_G(g)$ on $C_V(g)$. For this, put
\[\Irr_g(GV)\ =\ \{\chi\in\Irr(G)\ |\ \chi|_{\langle g\rangle\times C_V(g)}\mbox{ cannot be written as }
\rho_{\langle g\rangle}\times\psi\mbox{ for a character }\psi\mbox{ of }C_V(g)\}.\]
In particular, $\Irr(GV,g)\subseteq\Irr_g(GV)$.\\

\begin{nsatz}\label{satz3.1}
Let $G$ be a finite group and $V$ be a finite $G$--module. Let $g\in G$ such that $(o(g),|V|)=1$.
Write $A=\langle g\rangle$, $N=N_G(A)$ and
$C=C_V(g)$. Then  \\

(a) $|\Irr_g(GV)|\leq\frac{(n(N,A)-1)n(C_G(A),C)}{(|A|-1)|C|}\max_{1\not=a\in A}(|N_G(\langle
a\rangle)||C_V(a)|)$\\

(b) if $g$ is of prime order, then
\[|\Irr_g(GV)|\leq|C_G(A)|n(C_G(A),C)\]
(c) there are $X,Y\subseteq\Irr_g(GV)$ such that $\Irr_g(GV)$ is a disjoint union of $X$ and $Y$ and
\[|X|\leq\frac{(n(N,A)-1)n(C_G(A),C)}{(|A|-1)|C|^2}\max_{1\not=a\in A}(|N_G(\langle a\rangle)||C_V(a)|)
\mbox{ and}\]
\[|Y|\leq\frac{(n(N,A)-1)(n(C_G(A),C)-1)}{(|A|-1)|C|}\max_{1\not=a\in A}(|N_G(\langle a\rangle)||C_V(a)|)\]
(d) if $g$ is of prime order and $X,Y$ are as in (c), then
\[|X|\leq\frac{|C_G(A)|n(C_G(A),C)}{|C|}\mbox{ and }|Y|\leq|C_G(A)|(n(C_G(A),C)-1)\]
\end{nsatz}
\begin{beweis}
If $a_1,a_2\in A$ and $c_1,c_2\in C-\{1\}$, then it is straightforward to see that
$(a_1,c_1)^{GV}=(a_2,c_2)^{GV}$ implies that $a_1^G=a_2^G$. Hence if $T$ is a set of representatives
of the orbits of $N$ on $A-\{1\}$, then every conjugacy class of $GV$ that intersects nontrivially
with $(A-\{1\})\times C$ has a representative $ac$ for some $a\in T$ and some $c\in C$. Moreover, for each
$a\in T$ we have that if $c_3$, $c_4\in C$ are $C_G(A)$--conjugate, then $ac_3$ and $ac_4$ are
$C_G(A)$--conjugate
and thus $(ac_3)^G=(ac_4)^G$. This shows that for each $a\in T$ there are at most $n(C_G(A),C)$ conjugacy
classes of $GV$ intersecting nontrivially with $\{a\}\times C$.
Hence altogether we see that there are at most
\[|T|n(C_G(A),C)=(n(N,A)-1)n(C_G(A),C)\quad (1)\]
conjugacy classes of $GV$ which intersect $(A-\{1\})\times C$ nontrivially.\\
Moreover observe that for $1\not=a\in A$, $c\in C$, $h\in G$ and $u\in V$ we have
\[(ac)^{hu}\in A\times C\mbox{ if and only if }h\in N_G(\langle a\rangle),\ c^h\in C\mbox{ and }u\in C_V(a)\]
because the condition $(ac)^{hu}=a^h[a^h,u]c^h\in A\times C$ first forces $a^h\in A$ which implies
(as $A$ is cyclic) $a^h\in\langle a\rangle$, so $h\in N_G(\langle a\rangle)$, and then as
$c\in C\leq C_V(\langle a\rangle)$, it follows that $c^h\in C_V(\langle a\rangle)$ and
$[a^h,u]\in[\langle a\rangle, V]$. Now as by our hypothesis we have $V=C_V(\langle a\rangle)\times
[\langle a\rangle,V]$, we see that $(ac)^{hu}\in A\times C$ now forces $[a^h,u]=1$ and $c^h\in C$. Hence
$u\in C_V(a^h)=C_V(a)$.\\
Note that the direct product $A\times C$ is a subgroup of $GV$. We now define a generalized character
$\eta$ on $A\times C$ by
\[\eta=(|A|\cdot 1_A-\rho_A)\times 1_C\]
where $\rho_A$ is the regular character of $A$. So for $a\in A$, $c\in C$ we have
\[\eta(ac)=\left\{\begin{array}{ll}
0,&a=1\\
|A|,&a\not=1
\end{array}\right.\]
Therefore $\eta^{GV}$ vanishes on all conjugacy classes of $GV$ which intersect $(A-\{1\})\times C$ trivially,
whereas for $c\in C$ and $1\not=a\in A$ we have
\begin{eqnarray*}
\eta^{GV}(ac)&=&\frac{1}{|A\times C|}\sum_{{\scriptsize\begin{array}{l}
                                           h\in G\\
                                           u\in V
                                           \end{array}}}
                                           \dot{\eta}((ac)^{hu})\\
            &=&\frac{1}{|A||C|}\sum_{{\scriptsize\begin{array}{l}
                                    h\in N_G(\langle a\rangle)\\
                                    \mbox{with }c^h\in C
                                    \end{array}}}
                              \sum_{u\in C_V(a)}\eta((ac)^{hu})\\
           &=&\frac{1}{|A||C|}\sum_{{\scriptsize\begin{array}{l}
                                   h\in N_G(\langle a\rangle)\\
                                   \mbox{with }c^h\in C
                                   \end{array}}}
                              \sum_{u\in C_V(a)}\eta(a^hc^h)\\
           &=&\frac{|C_V(a)|}{|A||C|}\sum_{{\scriptsize\begin{array}{l}
                                          h\in N_G(\langle a\rangle)\\
                                          \mbox{with }c^h\in C
                                          \end{array}}}|A|\\
          &\leq&\frac{|N_G(\langle a\rangle)||C_V(a)|}{|C|}\qquad\qquad(2)
\end{eqnarray*}
Thus if $\{x_i\ |\ i=1,\ldots,k(GV)\}$ is a set of representatives for the conjugacy classes of $GV$,
then by (1) and (2) we see that
\begin{eqnarray*}
(n(N,A)-1)n(C_G(A),C)\cdot\frac{1}{|C|}\max_{1\not=a\in A}(|N_G(\langle a\rangle)||C_V(a)|)
      &\geq&\sum_{i=1}^{k(GV)}\eta^{GV}(x_i)\\
      &=&   \sum_{\tau\in\sIrr(GV)}(\tau\eta^{GV},\tau)_{GV}\\
      &=&   \sum_{\tau\in\sIrr(GV)}(\tau\eta,\tau)_{A\times C}\qquad (3).
\end{eqnarray*}
Observe that in case that $A$ is of prime order, then
\[n(N,A)-1=\frac{|A|-1}{|N:C_G(A)|}=\frac{(|A|-1)|C_G(A)|}{|N|}\]
and $\max\limits_{1\not=a\in A}(|N_G(\langle a\rangle)||C_V(a)|)=|N||C|$, so that (3) becomes
\[|C_G(A)|(|A|-1)n(C_G(A),C)\geq\sum_{\tau\in\sIrr(GV)}(\tau\eta,\tau)_{A\times C}\quad (3a)\]
Since $A\times C$ is a direct product, we can write
\[\tau_{A\times C}=\sum_{\lambda\in\sIrr(C)}(\tau_\lambda\times\lambda),\]
where $\tau_\lambda$ is a character of $A$ or $\tau_\lambda=0$. Then
\begin{eqnarray*}
(\tau\eta,\tau)_{A\times C}&=&\frac{1}{|A\times C|}\sum_{{\scriptsize\begin{array}{l}
                                                        a\in A\\
                                                        c\in C
                                                        \end{array}}}
                                                        \tau(ac)\eta(ac)\overline{\tau(ac)}\\
                       &=&\frac{1}{|A||C|}\sum_{{\scriptsize\begin{array}{l}
                                               1\not=a\in A\\
                                               c\in C
                                               \end{array}}}
                                               \tau(ac)|A|\overline{\tau(ac)}\\
                       &=&\frac{1}{|C|}\sum_{{\scriptsize\begin{array}{l}
                                            1\not=a\in A\\
                                            c\in C
                                            \end{array}}}\sum_{\lambda\in\sIrr(C)}
                                            \tau_\lambda(a)\lambda(c)\sum_{\mu\in\sIrr(C)}
                                            \overline{\tau_\mu (a)}\overline{\mu(c)}\\
                      &=&\sum_{1\not=a\in A}\sum_{\lambda,\mu\in\sIrr(C)}\tau_\lambda (a)\overline{\tau_\mu (a)}
                         \frac{1}{|C|}\sum_{c\in C}\lambda(c)\overline{\mu(c)}\\
                      &=&\sum_{1\not=a\in A}\sum_{\lambda,\mu\in\sIrr(C)}\tau_\lambda (a)\overline{\tau_\mu(a)}
                         (\lambda,\mu)_C
\end{eqnarray*}
As $(\lambda,\mu)_C=\left\{\begin{array}{ll}
                           1,&\lambda=\mu\\
                           0,&\lambda\not=\mu
                           \end{array}\right.$, we further obtain
\begin{eqnarray*}
(\tau\eta,\tau)_{A\times C}&=&\sum_{1\not=a\in A}\sum_{\lambda\in\sIrr(C)}\tau_\lambda (a)\overline{\tau_\lambda (a)}\\
                           &=&\sum_{\lambda\in\sIrr(C)}\sum_{1\not=a\in A}|\tau_\lambda(a)|^2\quad (4)
\end{eqnarray*}
Now observe that $\tau(1)=\sum\limits_{\lambda\in\Irr(C)}\tau_\lambda (1)$.\\
If all the $\tau_\lambda$ are multiples of $\rho_A$, then clearly $\tau_1\not\in\Irr_g(GV)$, and so if
$\tau\in\Irr_g(GV)$, then by \cite[Corollary 4]{robinson} with (4) we see that
\[(\tau\eta,\tau)_{A\times C}\geq |A|-1\quad (5)\]
So (3) and (5) yield
\[|\Irr_g(GV)|\leq\frac{(n(N,A)-1)n(C_G(A),C)}{(|A|-1)|C|}\max_{1\not=a\in A}(|N_G(\langle a\rangle)||C_V(a)|),
\quad (6)\]
and if $g$ is of prime order, then (3a) and (5) yield
\[|\Irr_g(GV)|\leq|C_G(A)|n(C_G(A),C).\quad (6a)\]
Now as in Section \ref{section2}, we now repeat the same arguments, but use
\[\eta_1=(|A| 1_A-\rho_A)\times (|C|1_C-\rho_C)\]
instead of $\eta$.\\
One can then easily check that
\[(n(N,A)-1)(n(C_G(A),C)-1)\cdot\frac{1}{|C|}\max_{1\not=a\in A}(|N_G(\langle a\rangle)||C_V(a)|)\geq
\sum_{\tau\in\sIrr(GV)}(\tau\eta_1,\tau)_{A\times C}\quad (3b)\]
and if $g$ is of prime order, then
\[|C_G(A)|(|A|-1)(n(C_G(A),C)-1)\geq\sum_{\tau\in\sIrr(GV)}(\tau\eta_1,\tau)_{A\times C}\quad (3c)\]
Moreover it is easily seen that
\begin{eqnarray*}
(\tau\eta_1,\tau)_{A\times C}&=&\sum_{{\scriptsize\begin{array}{l}
                                      1\not=a\in A\\
                                      1\not=c\in C
                                      \end{array}}}\tau(ac)\overline{\tau(ac)}\\
                            &=&\sum_{1\not=a\in A}\sum_{\lambda,\mu\in\sIrr(C)}\tau_\lambda (a)\overline{\tau_\mu (a)}
                               \sum_{1\not=c\in C}\lambda (c)\overline{\mu(c)},
\end{eqnarray*}
and as $\sum\limits_{1\not=c\in C}\lambda(c)\overline{\mu(c)}=\left\{\begin{array}{ll}
                                                                     -1,&\mbox{if }\lambda\not=\mu\\
                                                                     |C|-1,&\mbox{if }\lambda=\mu
                                                                     \end{array}\right. $,
it follows that
\[(\tau\eta_1,\tau)_{A\times C}=\sum_{\lambda<\mu}\sum_{1\not=a\in A}|\tau_\lambda (a)-\tau_\mu (a)|^2\quad (7)\]
where "$\leq$" is an arbitrary ordering on $\Irr(C)$.\\

Next suppose that there are exactly $a$ characters $\tau\in\Irr_g(GV)$ such that there is a character
$\psi$ of $A$ (depending on $\tau$) and there are $a_\lambda\in\z$ ($\lambda\in\Irr(C)$) such that
$\tau_\lambda=\psi+a_\lambda \rho_A$ for all $\lambda\in\Irr(C)$ and $\psi$ is not a multiple of $\rho_A$.
Then by (4) and \cite[Corollary 4]{robinson} we know that
\[(\tau\eta,\tau)_{A\times C}=\sum_{\lambda\in\sIrr(C)}\sum_{1\not=a\in A}|\psi(a)|^2\geq|C|(|A|-1)\]
and hence by (3) we get
\[a\leq\frac{(n(N,A)-1)n(C_G(A),C)}{(|A|-1)|C|^2}\max_{1\not=a\in A}(|N_G(\langle a\rangle )||C_V(a)|),\quad (8)\]
and if $g$ is of prime order, then by (3a) even
\[a\leq\frac{|C_G(A)|n(C_G(A),C)}{|C|}\quad (8a)\]
Now let $b$ be the number of $\tau\in\Irr_g(GV)$ such that there is no such $\psi$.\\
Then there exist $\lambda,\mu\in\Irr(C)$ with
\[\sum_{1\not=a\in A}|\tau_\lambda(a)-\tau_\mu(a)|^2\not=0,\]
and thus by \cite[Corollary 4]{robinson} we have
\[(\tau\eta_1,\tau)\geq |A|-1\quad (9)\]
So (3b) and (9) yield
\[b\leq\frac{(n(N,A)-1)(n(C_G(A),C)-1)}{|C|(|A|-1)}\max_{1\not=a\in A}(|N_G(\langle a\rangle)||C_V(a)|)\quad (10)\]
and, if $g$ is of prime order, then by (3c)
\[b\leq |C_G(A)|(n(C_G(A),C)-1),\quad (10b)\]
and clearly $a+b=|\Irr_g(GV)|$, and hence all the assertions follow and we are done.
\end{beweis}

\end{document}